\documentclass[reqno]{amsart}
\usepackage{amssymb}
\usepackage{hyperref}
\usepackage{lscape,graphicx}

\begin{document}
\title{Generalized Extension of Watson's theorem for the series $_{3}F_{2}(1)$}

\author[M. A. Rakha, M. M. Awad, A. O. Mohammed]
{Medhat A. Rakha, Mohammed M. Awad, Asmaa O. Mohammed}  % in alphabetical order

\address{Medhat A. Rakha \newline
 Department of Mathematics,
 Faculty of Science,
 Suez Canal University,
 El-Sheik Zayed 41522, 
 Ismailia - EGYPT}
\email{medhat\_rakha@scieence.suez.edu.eg}

\address{Mohammed M. Awad \newline
 Department of Mathematics,
 Faculty of Science,
 Suez Canal University,
 El-Sheik Zayed 41522, 
 Ismailia - EGYPT}
\email{mmawad\_sci@yahoo.com}

\address{Asmaa O. Mohammed \newline
 Department of Mathematics,
 Faculty of Science,
 Suez Canal University,
 El-Sheik Zayed 41522, 
 Ismailia - EGYPT}
\email{asmaa.orabi@science.suez.edu.eg}

%\thanks{Submitted December 8, 2011. Published March xx, 2012.}
\thanks{M. A. Rakha: Corresponding Author}


\begin{abstract}
The $_{3}F_{2}$ hypergeometric function plays a very significant role in the theory of hypergeometric and  generalized hypergeometric series. Despite that $_{3}F_{2}$ hypergeometric function has several applications in mathematics, also it has a lot of applications in physics and statistics. 

The fundamental purpose of this research paper is to find out the explicit expression of the $_{3}F_{2}$ Watson's classical summation theorem of the form:
\[
_{3}F_{2}\left[
\begin{array}
[c]{ccccc}%
a, & b, & c &  & \\
&  &  & ; & 1\\
\frac{1}{2}(a+b+i+1), & 2c+j &  &  &
\end{array}
\right]  
\]
with arbitrary $i$ and $j$, where for $i=j=0$, we get the well known Watson's theorem for the series $_{3}F_{2}(1)$.
\end{abstract}

\maketitle
\numberwithin{equation}{section}
\newtheorem{theorem}{Theorem}[section]
\newtheorem{lemma}[theorem]{Lemma}
\newtheorem{proposition}[theorem]{Proposition}
\newtheorem{corollary}[theorem]{Corollary}
\newtheorem*{remark}{Remark}
\section{Introduction}

The generalized hypergeometric function with $p$ numerator and $q$ denominator parameters is defined by \cite{1}

\begin{equation} 
_{p}F_{q}\left[
\begin{array}
[c]{ccccc}%
a_{1}, & ..., & a_{p} &  & \\
&  &  & ; & z\\
b_{1}, & ..., & b_{q} &  &
\end{array}
\right]  ={\displaystyle\sum\limits_{n=0}^{\infty}}
\frac{\left(  a_{1}\right)  _{n}...\left(  a_{p}\right)  _{n}}{\left(
b_{1}\right)  _{n}...\left(  b_{q}\right)  _{n}}\frac{z^{n}}{n!}, \label{1.1}
\end{equation} 

where $\left(a\right)_{n}$ denotes the shifted factorial defined for any
complex number $\alpha$, by
\[
\left(  \alpha\right)  _{n}=\left\{
\begin{array}
[c]{cc}
\alpha(\alpha+1)...(\alpha+n-1); & n=1,2,3,...\\
1; & n=0
\end{array}
\right.  .
\]

Using the main property $\Gamma\left(\alpha + 1\right)=\alpha \Gamma\left(\alpha \right)$, $\left(\alpha \right)_{n}$ can be written as
\[
\left(  \alpha\right)_{n}=\frac{\Gamma(\alpha+n)}{\Gamma(\alpha)}.
\]

It should by noted here that whenever hypergeometric and generalized hypergeometric functions summarized to be represented in term of Gamma function, the outcomes are critical from a theoretical and an appropriate point of view. Only some of summation theorems are available in literature, and it is Known as the classical summation theorems such as Gauss, Gauss's second, Kummer, and Bailey for the series $_{2}F_{1}$, Watson, Dixon, and Whipple for the series $_{3}F_{2}$.

The $_{3}F_{2}$ hypergeometric function plays a very remarkable part in the theory of hypergeometric and generalized hypergeometric series. Despite that $_{3}F_{2}$ hypergeometric function has several applications in mathematics such as in:
\begin{itemize}
\item[1.] Differential Equations: Since the generalized hypergeometric function and the monodromy of generalized hypergeometric function represented solutions of Picard-Fuchs equations, which used to solve many problems in classical mechanics and mathematical physics. For more information about such applications, see \cite{17}.

\item[2.] Conformal Mapping: Since in \cite{18} de Branges used the inequality 
\[
_{3}F_{2}\left[
\begin{array}
[c]{ccccc}%
$-n$, & n+\alpha, & \frac{1}{2}(\alpha+1) &  & \\
&  &  & ; & x\\
\alpha+1, & \frac{1}{2}(\alpha+3) &  &  &
\end{array}
\right] > 0, 
\]
where $ 0\leq x < 1,\alpha > -2$ and $n=0,1,2,...$, to proof the Bieberbach conjecture. Where the proof of this inequality is given in \cite{29}, see also \cite{30}.

\item[3.] Combinatorics and Number Theory: Many combinatorial identities are particular case of hypergeometric identities. For more details about such applications, see \cite{19},
\end{itemize}
also it has a lot of applications in physics and statistics such as: 
\begin{itemize}
\item[1.] Random Walks: Generalized hypergeometric function and Appell features appear in the assessment of the so-referred to as Watson integrals which symbolize the handiest possible lattice walks. They are also probably useful for the solution of extra complicated constrained lattice stroll issues. For further information about such application see \cite{20}.

\item[2.] Loop Integrals in Feynman Diagrams: Appell hypergeometric function gave One-loop integrals in Feynman diagrams also extension to two-loop, \cite{21, 22}.
 
\item[3.] $3j$, $6j$ and $9j$ Symbols: Since we can use $_{3}F_{2}$ functions with a unit argument to define the $3j$ symbols, which assume a critical part in the decay of reducible representations of the turning group into irreducible representations. Also recently, special cases of the $9j$ symbols are $_{5}F_{4}$ functions with a unit argument. Many  of combinatorial identities are individual cases of hypergeometric identities, see \cite{23, 24}.
\end{itemize}

Now, we begin by introducing the classical Watson's summation theorem $_{3}F_{2}$ of unit argument \cite{4}, which takes the form:

\begin{align}
&_{3}F_{2}\left[
\begin{array}
[c]{ccccc}%
a, & b, & c &  & \\
&  &  & ; & 1\\
\frac{1}{2}(a+b+1), & 2c &  &  &
\end{array}
\right] \nonumber\\
& =\frac{\Gamma(\frac{1}{2})\Gamma(c+\frac{1}{2})\Gamma(\frac{a}%
{2}+\frac{b}{2}+\frac{1}{2})\Gamma(c-\frac{a}{2}-\frac{b}{2}+\frac{1}{2}%
)}{\Gamma(\frac{a}{2}+\frac{1}{2})\Gamma(\frac{b}{2}+\frac{1}{2}%
)\Gamma(c-\frac{a}{2}+\frac{1}{2})\Gamma(c-\frac{b}{2}+\frac{1}{2})}\label{1.2}
\end{align}
where $Re(2c-a-b)>-1$.

 In \cite{2}, Watson gave the demonstrate of (\ref{1.2}) when one of the parameters $a$ or $b$ is a negative integer, and subsequently was established more generally in the non-terminating case by Whipple in \cite{3}.

The standard prove of (\ref{1.2}) given in \cite[p.149]{4} and \cite[p.54]{5}, depend on the following transformation due to Thomae  \cite{6}:
\begin{align*}
&_{3}F_{2}\left[
\begin{array}
[c]{ccccc}%
a, & b, & c &  & \\
&  &  & ; & 1\\
d, & e &  &  &
\end{array}
\right]\nonumber\\
  & =\frac{\Gamma(d)\Gamma(e)\Gamma(s)}{\Gamma(a)\Gamma(b+s)\Gamma
(c+s)} 
 ~\text{ }_{3}F_{2}\left[
\begin{array}
[c]{ccccc}%
d-a, & e-a, & s &  & \\
&  &  & ; & 1\\
b+s, & c+s &  &  &
\end{array}
\right]
\end{align*}
where $s=d+e-a-b-c$, $d=\frac{1}{2}(a+b+1)$ and $e=2c$.

MacRobert \cite{7} provided an alternative and more interested proof, by using the quadratic transformation for the Gauss's hypergeometric function \cite[Theorem 25, p. 67]{1}:
\begin{align*}
&_{2}F_{1}\left[
\begin{array}
[c]{cccc}%
2a, & 2b, &  & \\
&  & ; & x\\
a+b+\frac{1}{2} &  &  &
\end{array}
\right]=~_{2}F_{1}\left[
\begin{array}
[c]{cccc}%
a, & b, &  & \\
&  & ; & 4x(1-x)\\
a+b+\frac{1}{2} &  &  &
\end{array}
\right]
\end{align*}
valid for $\left\vert x\right\vert <1$ and $\left\vert 4x(1-x)\right\vert <1$.

Recently Rathie and Paris \cite{8} gave a basic confirmation of (\ref{1.2}) that just depends on the Gauss summation theorem for the $_{2}F_{1}$ hypergeometric function, namely, \cite{9}:

\[
_{2}F_{1}\left[
\begin{array}
[c]{cccc}%
a, & b, &  & \\
&  & ; & 1\\
c &  &  &
\end{array}
\right]  =\frac{\Gamma(c)\Gamma(c-a-b)}{\Gamma(c-a)\Gamma(c-b)}%
\]
where $Re(c-a-b)>0$, while Rakha in \cite{10} gave an extremely straightforward proof of(\ref{1.2}) by using the Gauss's second summation theorem:
\begin{align*}
&_{3}F_{2}\left[
\begin{array}
[c]{ccccc}%
a, & b, & c &  & \\
&  &  & ; & z\\
d, & 2b &  &  &
\end{array}
\right]\nonumber\\  &=\frac{\Gamma(d)}{\Gamma(c)\Gamma(d-s)}
\int_{0}^{1} t^{c-1}(1-t)^{d-c-1}~_{2}F_{1}\left[
\begin{array}
[c]{cccc}
a, & b, &  & \\
&  & ; & zt\\
2b &  &  &
\end{array}
\right]dt
\end{align*}

In 1987, Lavoie \cite{11}, obtained the following two summation formulas

\begin{align}
&_{3}F_{2}\left[
\begin{array}
[c]{ccccc}%
a, & b, & c &  & \\
&  &  & ; & 1\\
\frac{1}{2}(a+b+1), & 2c+1 &  &  &
\end{array}
\right] \nonumber\\  & =\frac{2^{a+b-2}\text{ }\Gamma(c+\frac{1}{2})\Gamma(\frac{a}%
{2}+\frac{b}{2}+\frac{1}{2})\Gamma(c-\frac{a}{2}-\frac{b}{2}+\frac{1}{2}%
)}{\Gamma(\frac{1}{2})\text{ }\Gamma(a)\text{ }\Gamma(b)} \nonumber \\
&  \times\left\{  \frac{\Gamma(\frac{a}{2})\Gamma(\frac{b}{2})}{\Gamma
(c-\frac{a}{2}+\frac{1}{2})\text{ }\Gamma(c-\frac{b}{2}+\frac{1}{2})}%
-\frac{\Gamma(\frac{a}{2}+\frac{1}{2})\text{ }\Gamma(\frac{b}{2}+\frac{1}{2}%
)}{\Gamma(c-\frac{a}{2}+1)\text{ }\Gamma(c-\frac{b}{2}+1)}\right\} \label{1.3}%
\end{align}
provided $Re(2c-a-b)>-3$, and

\begin{align}
&_{3}F_{2}\left[
\begin{array}
[c]{ccccc}%
a, & b, & c &  & \\
&  &  & ; & 1\\
\frac{1}{2}(a+b+1), & 2c-1 &  &  &
\end{array}
\right]\nonumber\\  &=\frac{2^{a+b-2}\text{ }\Gamma(c-\frac{1}{2})\Gamma(\frac{a}%
{2}+\frac{b}{2}+\frac{1}{2})\Gamma(c-\frac{a}{2}-\frac{b}{2}-\frac{1}{2}%
)}{\Gamma(\frac{1}{2})\text{ }\Gamma(a)\text{ }\Gamma(b)}\nonumber \\
&  \times\left\{  \frac{\Gamma(\frac{a}{2})\Gamma(\frac{b}{2})}{\Gamma
(c-\frac{a}{2}-\frac{1}{2})\text{ }\Gamma(c-\frac{b}{2}-\frac{1}{2})}%
-\frac{\Gamma(\frac{a}{2}+\frac{1}{2})\text{ }\Gamma(\frac{b}{2}+\frac{1}{2}%
)}{\Gamma(c-\frac{a}{2})\text{ }\Gamma(c-\frac{b}{2})}\right\} 
\label{1.4}%
\end{align}
provided $Re(2c-a-b)>1.$\\

In 1992, Lavoie et al. in \cite{12}, took out explicit expression of the series

\begin{equation}
_{3}F_{2}\left[
\begin{array}
[c]{ccccc}%
a, & b, & c &  & \\
&  &  & ; & 1\\
\frac{1}{2}(a+b+i+1), & 2c+j &  &  &
\end{array}
\right] \label{1.5}%
\end{equation}
for $i,j=0,\pm1,\pm2$, where at $i=j=0$ we obtain (\ref{1.2}).

Other remarkable results of such computations, are:
%\begin{description}
\begin{itemize}
\item[1.] Stainislaw \cite{13}, in 1997 gave an analytical formula for (\ref{1.5}) with a fixed $j$ and arbitrary $i$.

\item[2.]  Kim et al., in \cite{14}, have obtained the above result (\ref{1.5}) for $j=0$ and $i=0,\pm1,\pm2,...,\pm5$.

\item[3.] Chu \cite{15}, in 2012, investigates the generalized Watson's series with two extra integer parameters by combining the linearization method with Dougall's sum for well-poised $_{5}F_{4}$-series.

\item[4.] in 2013, Rakha et al. in \cite{16}, established result (\ref{1.5}) for $i=0,\pm1,\pm2,...,\pm5;$ $j=0,\pm 1,\pm 2$.
\end{itemize}
%\end{description}

The major purpose of this paper is to find explicit extensions of the classical Watson's summation theorem (\ref{1.5}) for any $i$ and $j$, to have more summations theorems as well as more contiguous relations about the $_{3}F_{2}(1)$, hypergeometric series.
 
\section{Main Results}
Let us consider that
\begin{align}
f_{i,j}(a,b,c)  &=~_{3}F_{2}\left[
\begin{array}
[c]{ccccc}%
a & b & c &  & \\
&  &  & ; & 1\\
\frac{a+b+i+1}{2} & 2c+j &  &  &
\end{array}
\right]\nonumber\\ 
&= {\displaystyle\sum\limits_{n=0}^{\infty}}
\frac{\left(  a\right)  _{n}\left(  b\right)  _{n}\left(  c\right)  _{n}
}{\left(  \frac{a+b+i+1}{2}\right)  _{n}\left(  2c+j\right)  _{n}}\frac{1}{n!}.\label{1}
\end{align}

It is clear that
\begin{align*}
&(2c+j)\text{ }f_{i,j+1}(a,b,c)\\
& ={\displaystyle\sum\limits_{n=0}^{\infty}}
\frac{\left(  a\right)  _{n}\left(  b\right)  _{n}\left(  c\right)
_{n}(2c+j)}{\left(  \frac{a+b+i+1}{2}\right)  _{n}\left(  2c+j+1\right)  _{n}}\frac{1}{n!}\\
& ={\displaystyle\sum\limits_{n=0}^{\infty}}
\frac{\left(  a\right)  _{n}\left(  b\right)  _{n}\left(  c\right)
_{n}(2c+j+n)}{\left(  \frac{a+b+i+1}{2}\right)  _{n}\left(  2c+j+1\right)
_{n}}\frac{1}{n!}-{\displaystyle\sum\limits_{n=0}^{\infty}}
\frac{\left(  a\right)  _{n}\left(  b\right)  _{n}\left(  c\right)  _{n}\text{}n}{\left(  \frac{a+b+i+1}{2}\right)  _{n}\left(  2c+j+1\right)  _{n}}\frac{1}{n!}\\
& =(2c+j){\displaystyle\sum\limits_{n=0}^{\infty}}
\frac{\left(  a\right)  _{n}\left(  b\right)  _{n}\left(  c\right)  _{n}%
}{\left(  \frac{a+b+i+1}{2}\right)  _{n}\left(  2c+j\right)  _{n}}\frac{1}{n!}-{\displaystyle\sum\limits_{n=0}^{\infty}}
\frac{\left(  a\right)  _{n+1}\left(  b\right)  _{n+1}\left(  c\right)
_{n+1}\text{ }}{\left(  \frac{a+b+i+1}{2}\right)  _{n+1}\left(  2c+j+1\right)_{n+1}}\frac{1}{n!}\\
& =(2c+j)\text{ }f_{i,j}(a,b,c)-\frac{2abc}{(a+b+i+1)(2c+j+1)}
{\displaystyle\sum\limits_{n=0}^{\infty}}
\frac{\left(  a+1\right)  _{n+1}\left(  b+1\right)  _{n+1}\left(  c+1\right)_{n+1}\text{ }}{\left(  \frac{a+b+i+3}{2}\right)  _{n+1}\left(  2c+j+2\right)_{n+1}}\frac{1}{n!}\\
& =(2c+j)\text{ }f_{i,j}(a,b,c)-\frac{2abc}{(a+b+i+1)(2c+j+1)}f_{i,j}%
(a+1,b+1,c+1),
\end{align*}
from which, we conclude the following general important main relation:
\begin{align}
& (2c+j)\text{ }f_{i,j+1}(a,b,c) \nonumber\\
& =(2c+j)\text{ }f_{i,j}(a,b,c)-\frac
{2abc}{(a+b+i+1)(2c+j+1)}f_{i,j}(a+1,b+1,c+1).\label{1.6}%
\end{align}

So, if we know $f_{i,0}(a,b,c)$, we can generate $f_{i,j}(a,b,c)$ for any values of $i$ and $j$.

\subsection{Special Cases}
%\begin{itemize}
\begin{enumerate}
\item When $i=j=0$ in (\ref{1.6}), we obtain
\begin{align*}
&_{3}F_{2}\left[
\begin{array}
[c]{ccccc}%
a, & b, & c &  & \\
&  &  & ; & 1\\
\frac{1}{2}(a+b+1), & 2c+1 &  &  &
\end{array}
\right] \nonumber\\
 & =~_{3}F_{2}\left[
\begin{array}
[c]{ccccc}%
a, & b, & c &  & \\
&  &  & ; & 1\\
\frac{1}{2}(a+b+1), & 2c &  &  &
\end{array}
\right] \nonumber\\
&  -\frac{ab}{(2c+1)(a+b+1)}\,_{3}F_{2}\left[
\begin{array}
[c]{ccccc}%
a+1, & b+1, & c+1 &  & \\
&  &  & ; & 1\\
\frac{1}{2}(a+b+3), & 2c+2 &  &  &
\end{array}
\right]
\end{align*}
which appeared in \cite[Eq.(3.18), pp. 229]{16}, \cite[Result (2), pp.269]{11},\cite[Result (1),pp.24]{12}, \cite[Eq.(4.4),pp.12]{25} and \cite[Theorem 4, p.147]{28}.
\begin{enumerate}
\item In such a case, the result when $a=1$, $b=1$ and $c=1$, appeared in \cite[Result 209, p. 459]{26}.
\item In such a case the result when $a=\frac{2}{3}$, $b=\frac{4}{3}$ and $c=1$, appeared in \cite[Eq. 26]{27}.
\end{enumerate}
%\end{itemize}

\item When $i=1$ and $j=0$ in (\ref{1.6}), we obtain
\begin{align*}
&_{3}F_{2}\left[
\begin{array}
[c]{ccccc}%
a, & b, & c &  & \\
&  &  & ; & 1\\
\frac{1}{2}(a+b+2), & 2c+1 &  &  &
\end{array}
\right] \nonumber\\
 & =~_{3}F_{2}\left[
\begin{array}
[c]{ccccc}%
a, & b, & c &  & \\
&  &  & ; & 1\\
\frac{1}{2}(a+b+2), & 2c &  &  &
\end{array}
\right] \nonumber\\
&  -\frac{ab}{(2c+1)(a+b+2)}\,_{3}F_{2}\left[
\begin{array}
[c]{ccccc}%
a+1, & b+1, & c+1 &  & \\
&  &  & ; & 1\\
\frac{1}{2}(a+b+4), & 2c+2 &  &  &
\end{array}
\right]
\end{align*}
which appeared in \cite[Eq.(3.18), pp.229]{16}, \cite[Theorem 5, pp.148]{28}, \cite[Result (1), pp.24]{12}, and \cite[Theorem 2, pp. 144]{28}.

In such a case, the results when  $a=\frac{1}{2}$, $b=1$, $c=\frac{5}{4}$; $a=\frac{1}{2}$, $b=\frac{3}{2}$, $c=1$ and $a=b=c=1$;  appeared in \cite[Results 187, 188 \& 211, pages 458, 458 \& 459]{26}, respectively.

\item When $i=2$ and $j=0$ in (\ref{1.6}), we obtain
\begin{align*}
&_{3}F_{2}\left[
\begin{array}
[c]{ccccc}%
a, & b, & c &  & \\
&  &  & ; & 1\\
\frac{1}{2}(a+b+3), & 2c+1 &  &  &
\end{array}
\right] \nonumber\\
 & =~_{3}F_{2}\left[
\begin{array}
[c]{ccccc}%
a, & b, & c &  & \\
&  &  & ; & 1\\
\frac{1}{2}(a+b+3), & 2c &  &  &
\end{array}
\right] \nonumber\\
&  -\frac{ab}{(2c+1)(a+b+3)}\,_{3}F_{2}\left[
\begin{array}
[c]{ccccc}%
a+1, & b+1, & c+1 &  & \\
&  &  & ; & 1\\
\frac{1}{2}(a+b+5), & 2c+2 &  &  &
\end{array}
\right]
\end{align*}
which appeared in \cite[Eq.3.18, pp.229]{16},\cite[Eq.(4.5),pp.12]{25} and \cite[Result(1),pp.24]{12}.

In such a case, the results when $a=1$, $b=\frac{3}{2}$, $c=\frac{3}{4}$;  $a=1$, $b=2$, $c=1$ and $a=\frac{3}{2}$, $b=\frac{3}{2}$, $c=1$; appeared in \cite[Results 204, 234 \& 242, pages 459 \& 460]{26}, respectively.

\item When $i=3$ and $j=0$ in (\ref{1.6}), we obtain
\begin{align*}
&_{3}F_{2}\left[
\begin{array}
[c]{ccccc}%
a, & b, & c &  & \\
&  &  & ; & 1\\
\frac{1}{2}(a+b+4), & 2c+1 &  &  &
\end{array}
\right] \nonumber\\
 & =~_{3}F_{2}\left[
\begin{array}
[c]{ccccc}%
a, & b, & c &  & \\
&  &  & ; & 1\\
\frac{1}{2}(a+b+4), & 2c &  &  &
\end{array}
\right] \nonumber\\
&  -\frac{ab}{(2c+1)(a+b+4)}\,_{3}F_{2}\left[
\begin{array}
[c]{ccccc}%
a+1, & b+1, & c+1 &  & \\
&  &  & ; & 1\\
\frac{1}{2}(a+b+6), & 2c+2 &  &  &
\end{array}
\right]
\end{align*}
which appeared in \cite[Eq.(3.18), pp.229]{16} and \cite[Result(2.23), pp.380]{13}.

In such a case, the result when $a=1$, $b=3$ and $c=1$, appeared in \cite[Result (237), p.460]{26}.

\item When $i=4$ and $j=0$ in (\ref{1.6}), we obtain
\begin{align*}
&_{3}F_{2}\left[
\begin{array}
[c]{ccccc}%
a, & b, & c &  & \\
&  &  & ; & 1\\
\frac{1}{2}(a+b+5), & 2c+1 &  &  &
\end{array}
\right] \nonumber\\
 & =~_{3}F_{2}\left[
\begin{array}
[c]{ccccc}%
a, & b, & c &  & \\
&  &  & ; & 1\\
\frac{1}{2}(a+b+5), & 2c &  &  &
\end{array}
\right] \nonumber\\
&  -\frac{ab}{(2c+1)(a+b+5)}\,_{3}F_{2}\left[
\begin{array}
[c]{ccccc}%
a+1, & b+1, & c+1 &  & \\
&  &  & ; & 1\\
\frac{1}{2}(a+b+7), & 2c+2 &  &  &
\end{array}
\right]
\end{align*}
which appeared in \cite[Eq. (3.18), pp.229]{16}.

In such a case, the result when $a=4$, $b=c=1$ appeared in  \cite[Results 240 \& 239, page 460]{26}.

 \item When $i=5$ and $j=0$ in (\ref{1.6}), we obtain
\begin{align*}
&_{3}F_{2}\left[
\begin{array}
[c]{ccccc}%
a, & b, & c &  & \\
&  &  & ; & 1\\
\frac{1}{2}(a+b+6), & 2c+1 &  &  &
\end{array}
\right] \nonumber\\
 & =~_{3}F_{2}\left[
\begin{array}
[c]{ccccc}%
a, & b, & c &  & \\
&  &  & ; & 1\\
\frac{1}{2}(a+b+6), & 2c &  &  &
\end{array}
\right] \nonumber\\
&  -\frac{ab}{(2c+1)(a+b+6)}\,_{3}F_{2}\left[
\begin{array}
[c]{ccccc}%
a+1, & b+1, & c+1 &  & \\
&  &  & ; & 1\\
\frac{1}{2}(a+b+8), & 2c+2 &  &  &
\end{array}
\right]
\end{align*}
which appeared in \cite[Eq.(3.18), pp.229]{16} and \cite[Result (2.22), pp.380]{13}.

In such a case, the results when $a=c=1$, $b=3$ and $a=c=1$, $b=5$ appeared in \cite[Results 238 \& 241, p. 460]{26}, respectively.

\item When $i=-1$ and $j=0$ in (\ref{1.6}), we obtain
\begin{align*}
&_{3}F_{2}\left[
\begin{array}
[c]{ccccc}%
a, & b, & c &  & \\
&  &  & ; & 1\\
\frac{1}{2}(a+b), & 2c+1 &  &  &
\end{array}
\right] \nonumber\\
 & =~_{3}F_{2}\left[
\begin{array}
[c]{ccccc}%
a, & b, & c &  & \\
&  &  & ; & 1\\
\frac{1}{2}(a+b), & 2c &  &  &
\end{array}
\right] \nonumber\\
&  -\frac{ab}{(2c+1)(a+b)}\,_{3}F_{2}\left[
\begin{array}
[c]{ccccc}%
a+1, & b+1, & c+1 &  & \\
&  &  & ; & 1\\
\frac{1}{2}(a+b+2), & 2c+2 &  &  &
\end{array}
\right]
\end{align*}
which appeared in \cite[Result (1), pp.24]{12}, \cite[Eq.(3.18), pp.229]{16}, \cite[Example (11), pp.9]{15} and \cite[Theorem (1), pp.143]{28}.

\item When $i=-2$ and $j=0$ in (\ref{1.6}), we obtain
\begin{align*}
&_{3}F_{2}\left[
\begin{array}
[c]{ccccc}%
a, & b, & c &  & \\
&  &  & ; & 1\\
\frac{1}{2}(a+b-1), & 2c+1 &  &  &
\end{array}
\right] \nonumber\\
 & =~_{3}F_{2}\left[
\begin{array}
[c]{ccccc}%
a, & b, & c &  & \\
&  &  & ; & 1\\
\frac{1}{2}(a+b-1), & 2c &  &  &
\end{array}
\right] \nonumber\\
&  -\frac{ab}{(2c+1)(a+b-1)}\,_{3}F_{2}\left[
\begin{array}
[c]{ccccc}%
a+1, & b+1, & c+1 &  & \\
&  &  & ; & 1\\
\frac{1}{2}(a+b+1), & 2c+2 &  &  &
\end{array}
\right]
\end{align*}
which appeared in \cite[Result(1), pp.24]{12} and \cite[Eq. (3.18), pp.229]{16}.

\item When $i=-3$ and $j=0$ in (\ref{1.6}), we obtain
\begin{align*}
&_{3}F_{2}\left[
\begin{array}
[c]{ccccc}%
a, & b, & c &  & \\
&  &  & ; & 1\\
\frac{1}{2}(a+b-2), & 2c+1 &  &  &
\end{array}
\right] \nonumber\\
 & =~_{3}F_{2}\left[
\begin{array}
[c]{ccccc}%
a, & b, & c &  & \\
&  &  & ; & 1\\
\frac{1}{2}(a+b-2), & 2c &  &  &
\end{array}
\right] \nonumber\\
&  -\frac{ab}{(2c+1)(a+b-2)}\,_{3}F_{2}\left[
\begin{array}
[c]{ccccc}%
a+1, & b+1, & c+1 &  & \\
&  &  & ; & 1\\
\frac{1}{2}(a+b), & 2c+2 &  &  &
\end{array}
\right]
\end{align*}
which appeared in \cite[Eq.(3.18), pp.229]{16}.

\item When $i=-4$ and $j=0$ in (\ref{1.6}), we obtain
\begin{align*}
&_{3}F_{2}\left[
\begin{array}
[c]{ccccc}%
a, & b, & c &  & \\
&  &  & ; & 1\\
\frac{1}{2}(a+b-3), & 2c+1 &  &  &
\end{array}
\right] \nonumber\\
 & =~_{3}F_{2}\left[
\begin{array}
[c]{ccccc}%
a, & b, & c &  & \\
&  &  & ; & 1\\
\frac{1}{2}(a+b-3), & 2c &  &  &
\end{array}
\right] \nonumber\\
&  -\frac{ab}{(2c+1)(a+b-3)}\,_{3}F_{2}\left[
\begin{array}
[c]{ccccc}%
a+1, & b+1, & c+1 &  & \\
&  &  & ; & 1\\
\frac{1}{2}(a+b-1), & 2c+2 &  &  &
\end{array}
\right]
\end{align*}
which appeared in \cite[Eq.(3.18), pp.229]{16}.

\item When $i=-5$ and $j=0$ in (\ref{1.6}), we obtain
\begin{align*}
&_{3}F_{2}\left[
\begin{array}
[c]{ccccc}%
a, & b, & c &  & \\
&  &  & ; & 1\\
\frac{1}{2}(a+b-4), & 2c+1 &  &  &
\end{array}
\right] \nonumber\\
 & =~_{3}F_{2}\left[
\begin{array}
[c]{ccccc}%
a, & b, & c &  & \\
&  &  & ; & 1\\
\frac{1}{2}(a+b-4), & 2c &  &  &
\end{array}
\right] \nonumber\\
&  -\frac{ab}{(2c+1)(a+b-4)}\,_{3}F_{2}\left[
\begin{array}
[c]{ccccc}%
a+1, & b+1, & c+1 &  & \\
&  &  & ; & 1\\
\frac{1}{2}(a+b-2), & 2c+2 &  &  &
\end{array}
\right]
\end{align*}
which appeared in \cite[Eq.(3.18), pp.229]{16}.

\item When $i=0$ and $j=1$ in (\ref{1.6}), we obtain
\begin{align*}
&_{3}F_{2}\left[
\begin{array}
[c]{ccccc}%
a, & b, & c &  & \\
&  &  & ; & 1\\
\frac{1}{2}(a+b+1), & 2c+2 &  &  &
\end{array}
\right] \nonumber\\
 & =~_{3}F_{2}\left[
\begin{array}
[c]{ccccc}%
a, & b, & c &  & \\
&  &  & ; & 1\\
\frac{1}{2}(a+b+1), & 2c+1 &  &  &
\end{array}
\right] \nonumber\\
&  -2\frac{abc}{(2c+1)(2c+2)(a+b+1)}\,_{3}F_{2}\left[
\begin{array}
[c]{ccccc}%
a+1, & b+1, & c+1 &  & \\
&  &  & ; & 1\\
\frac{1}{2}(a+b+3), & 2c+3 &  &  &
\end{array}
\right]
\end{align*}
which appeared in \cite[Eq.(3.20), pp. 230]{16} and \cite[Result (1), pp.24]{12}.

In such a case, the results when $a=\frac{1}{4}$, $b=1$, $c=\frac{1}{8}$;$a=\frac{1}{3}$, $b=1$, $c=\frac{1}{6}$; $a=\frac{1}{2}$, $b=1$ $c=\frac{1}{4}$; $a=\frac{1}{2}$, $b=1$ $c=\frac{3}{8}$; $a=\frac{3}{4}$, $b=1$ $c=\frac{3}{8}$ and $a=\frac{3}{2}$, $b=1$ $c=\frac{3}{4}$ appeared in \cite[Results 123, 130, 136, 137, 160 \& 203, pages 456, 457 \& 459]{26}.

\item When $i=1$ and $j=1$ in (\ref{1.6}) we obtain
\begin{align*}
&_{3}F_{2}\left[
\begin{array}
[c]{ccccc}%
a, & b, & c &  & \\
&  &  & ; & 1\\
\frac{1}{2}(a+b+2), & 2c+2 &  &  &
\end{array}
\right] \nonumber\\
 & =~_{3}F_{2}\left[
\begin{array}
[c]{ccccc}%
a, & b, & c &  & \\
&  &  & ; & 1\\
\frac{1}{2}(a+b+2), & 2c+1 &  &  &
\end{array}
\right] \nonumber\\
&  -2\frac{abc}{(2c+1)(2c+2)(a+b+2)}\,_{3}F_{2}\left[
\begin{array}
[c]{ccccc}%
a+1, & b+1, & c+1 &  & \\
&  &  & ; & 1\\
\frac{1}{2}(a+b+4), & 2c+3 &  &  &
\end{array}
\right]
\end{align*}
which appeared in \cite[Result (1), pp. 24]{12}, \cite[Theorem (5), pp.148]{28} and \cite[Example (9), pp.8]{15}.

In such a case, the result when $a=\frac{3}{2}$, $b=1$ and $c=\frac{1}{4}$, appeared in \cite[Result 148,p.457]{26}.

\item When $i=2$ and $j=1$ in (\ref{1.6}), we obtain
\begin{align*}
&_{3}F_{2}\left[
\begin{array}
[c]{ccccc}%
a, & b, & c &  & \\
&  &  & ; & 1\\
\frac{1}{2}(a+b+3), & 2c+2 &  &  &
\end{array}
\right] \nonumber\\
 & =~_{3}F_{2}\left[
\begin{array}
[c]{ccccc}%
a, & b, & c &  & \\
&  &  & ; & 1\\
\frac{1}{2}(a+b+3), & 2c+1 &  &  &
\end{array}
\right] \nonumber\\
&  -2\frac{abc}{(2c+1)(2c+2)(a+b+3)}\,_{3}F_{2}\left[
\begin{array}
[c]{ccccc}%
a+1, & b+1, & c+1 &  & \\
&  &  & ; & 1\\
\frac{1}{2}(a+b+5), & 2c+3 &  &  &
\end{array}
\right]
\end{align*}
which appeared in \cite[Result(1), pp.24]{12}.

In such a case, the result when $a=\frac{1}{2}$, $b=1$ and $c=\frac{1}{4}$, appeared in \cite[Result 139, p. 456]{26}.

\item When $i=0$ and $j=-1$ in (\ref{1.6}), we obtain
\begin{align*}
&_{3}F_{2}\left[
\begin{array}
[c]{ccccc}%
a, & b, & c &  & \\
&  &  & ; & 1\\
\frac{1}{2}(a+b+1), & 2c-1 &  &  &
\end{array}
\right] \nonumber\\
 & =~_{3}F_{2}\left[
\begin{array}
[c]{ccccc}%
a, & b, & c &  & \\
&  &  & ; & 1\\
\frac{1}{2}(a+b+1), & 2c &  &  &
\end{array}
\right] \nonumber\\
&  +\frac{ab}{(2c-1)(a+b+1)}\,_{3}F_{2}\left[
\begin{array}
[c]{ccccc}%
a+1, & b+1, & c+1 &  & \\
&  &  & ; & 1\\
\frac{1}{2}(a+b+3), & 2c+1 &  &  &
\end{array}
\right]
\end{align*}
which appeared in \cite[Result(1),pp.24]{12}, \cite[Eq.(3.19), pp. 230]{16}, \cite[Result(1), pp. 269]{11} and \cite[Theorem (7), pp. 152]{28}.

In such a case, the result when $a=\frac{1}{2}$, $b=\frac{1}{2}$ and $c=2$, appeared in \cite[Result 172, p.458]{26} \& \cite[Result(2.9),pp.5]{31}. 

 \item When $i=1$ and $j=-1$ in (\ref{1.6}), we obtain
\begin{align*}
&_{3}F_{2}\left[
\begin{array}
[c]{ccccc}%
a, & b, & c &  & \\
&  &  & ; & 1\\
\frac{1}{2}(a+b+2), & 2c-1 &  &  &
\end{array}
\right] \nonumber\\
 & =~_{3}F_{2}\left[
\begin{array}
[c]{ccccc}%
a, & b, & c &  & \\
&  &  & ; & 1\\
\frac{1}{2}(a+b+2), & 2c &  &  &
\end{array}
\right] \nonumber\\
&  +\frac{ab}{(2c-1)(a+b+2)}\,_{3}F_{2}\left[
\begin{array}
[c]{ccccc}%
a+1, & b+1, & c+1 &  & \\
&  &  & ; & 1\\
\frac{1}{2}(a+b+4), & 2c+1 &  &  &
\end{array}
\right]
\end{align*}
which appeared in \cite[Result(1), pp. 24]{12},\cite[Eq.(3.20), pp. 230]{16}, \cite[Example(10), pp. 8]{15} and \cite[Theorem (8), pp.153]{28}.

\item When $i=2$ and $j=-1$ in (\ref{1.6}), we obtain
\begin{align*}
&_{3}F_{2}\left[
\begin{array}
[c]{ccccc}%
a, & b, & c &  & \\
&  &  & ; & 1\\
\frac{1}{2}(a+b+3), & 2c-1 &  &  &
\end{array}
\right] \nonumber\\
 & =~_{3}F_{2}\left[
\begin{array}
[c]{ccccc}%
a, & b, & c &  & \\
&  &  & ; & 1\\
\frac{1}{2}(a+b+3), & 2c &  &  &
\end{array}
\right] \nonumber\\
&  +\frac{ab}{(2c-1)(a+b+3)}\,_{3}F_{2}\left[
\begin{array}
[c]{ccccc}%
a+1, & b+1, & c+1 &  & \\
&  &  & ; & 1\\
\frac{1}{2}(a+b+4), & 2c+1 &  &  &
\end{array}
\right]
\end{align*}
which appeared in \cite[Result(1), pp.24]{12}.

In such a case, the results when $a=1$, $b=2$ and $c=2$  appeared in \cite[Results 243, page 460]{26}. 

\item When $i=1$ and $j=-2$ in (\ref{1.6}), we obtain
\begin{align*}
&_{3}F_{2}\left[
\begin{array}
[c]{ccccc}%
a, & b, & c &  & \\
&  &  & ; & 1\\
\frac{1}{2}(a+b+2), & 2c-1 &  &  &
\end{array}
\right] \nonumber\\
 & =~_{3}F_{2}\left[
\begin{array}
[c]{ccccc}%
a, & b, & c &  & \\
&  &  & ; & 1\\
\frac{1}{2}(a+b+2), & 2c-2 &  &  &
\end{array}
\right] \nonumber\\
&  -\frac{abc}{(2c-1)(c-1)(a+b+2)}\,_{3}F_{2}\left[
\begin{array}
[c]{ccccc}%
a+1, & b+1, & c+1 &  & \\
&  &  & ; & 1\\
\frac{1}{2}(a+b+4), & 2c &  &  &
\end{array}
\right]
\end{align*}
which appeared in \cite[Result(1), pp.24]{12}, \cite[Example(10), pp. 8]{15} and \cite[Theorem(8), pp. 153]{28}.

\item When $i=-1$ and $j=-2$ in (\ref{1.6}), we obtain
\begin{align*}
&_{3}F_{2}\left[
\begin{array}
[c]{ccccc}%
a, & b, & c &  & \\
&  &  & ; & 1\\
\frac{1}{2}(a+b), & 2c-1 &  &  &
\end{array}
\right] \nonumber\\
 & =~_{3}F_{2}\left[
\begin{array}
[c]{ccccc}%
a, & b, & c &  & \\
&  &  & ; & 1\\
\frac{1}{2}(a+b), & 2c-2 &  &  &
\end{array}
\right] \nonumber\\
&  -\frac{abc}{(2c-1)(c-1)(a+b)}\,_{3}F_{2}\left[
\begin{array}
[c]{ccccc}%
a+1, & b+1, & c+1 &  & \\
&  &  & ; & 1\\
\frac{1}{2}(a+b+2), & 2c &  &  &
\end{array}
\right]
\end{align*}
which appeared in \cite[Result(1), pp.24]{12}, \cite[Example(12), pp. 9]{15} and \cite[Theorem(6), pp. 150]{28}.
\end{enumerate}

%\end{itemize} 

\section{Concluding Remarks}
\begin{itemize}
\item[1.] Various other special cases of our result can be obtained.
\item[2.] Many new identities and relations which obtained from our result are under examinations and will be published later.
\item[3.] We have already established in the previous section a recursive relation (\ref{1.6}), that generalized the extension of Watson summation theorem $_{3}F_{2}(1)$. Another explicit expression of (\ref{1.5}) that generalize our result (\ref{1.6}), can be presented in the next theorem.
\begin{theorem}
 \begin{align*}
f_{i,j}(a,b,c) & =(2c+j)\text{ }f_{i+1,j}(a-1,b,c)\\
&-\frac{2ab}{(a+b+i+1)(2c+j)}f_{i+1,j-1}(a,b+1,c+1)
\end{align*}
\end{theorem}
where $f_{i,j}(a,b,c)$ is defined as in (\ref{1}).
\end{itemize}

\subsection*{Conflict of Interests}

The authors declare that they have no any conflict of interests.

\subsection*{Acknowledgments}
All authors contributed equally in this paper. They read and approved the final manuscript.

\end{document}